\input amstex
\documentstyle{amsppt}
\magnification=\magstep1
\pageheight{9truein}
\pagewidth{6.7truein}
\NoBlackBoxes
%%%%%%%%%%%%%%%%%%%%%%%%%%%%%%%%%%
\def\ds{\displaystyle}
\def\ep{\varepsilon} 
 
\def\A{{\Cal A}}

\def\nat{{\Bbb N}}
\def\que{{\Bbb Q}}
\def\real{{\Bbb R}}
\def\normm{{|\!|\!|\,}}  

\def\dotNorm{\|\cdot\|}
\def\dotnormm{\normm\cdot\normm}
\def\supp{\operatorname{supp}}
%%%%%%%%%%%%%%%%%%%%%%%%%%%%%%%%%%

\topmatter
\title Asymptotic Properties of Banach Spaces under Renormings \endtitle
\rightheadtext{Asymptotic Properties under Renormings}

\author E. Odell and Th. Schlumprecht\endauthor
\address Department of Mathematics, The University of Texas at Austin, 
Austin, TX 78712-1082\endaddress
\email odell\@math.utexas.edu\endemail
\address Department of Mathematics, Texas A\&M University, College 
Station, TX 77843-3368\endaddress
\email schlump\@math.tamu.edu\endemail 
\thanks Research of both authors was supported by NSF and TARP.\endthanks 
\abstract 
It is shown that a separable Banach space $X$ can be given an equivalent norm 
$\dotnormm$ with the following properties:\quad 
If $(x_n)\subseteq X$ is relatively 
weakly compact and $\lim_{m\to\infty} \lim_{n\to\infty}\break 
\normm x_m+x_n\normm 
= 2\lim_{m\to\infty} \normm x_m\normm$ then $(x_n)$ converges in norm. 
This yields a characterization of reflexivity once proposed by V.D.~Milman.
In addition it is shown that some spreading model of a 
sequence in $(X,\dotnormm)$ is 1-equivalent to the unit vector basis 
of $\ell_1$ (respectively, $c_0$) implies that 
$X$ contains an isomorph of $\ell_1$ 
(respectively, $c_0$).
\endabstract 
\keywords spreading model, Ramsey theory, $\ell_1$, $c_0$, reflexive 
Banach space\endkeywords
\subjclass 46B03, 46B45\endsubjclass
\endtopmatter

\document 
\baselineskip=18pt 

\head \S1. Introduction\endhead 

A classical problem in functional analysis has been to give a geometric 
characterization of reflexivity for a Banach space. 
The first result of this type was D.P.~Milman's \cite{Mil} 
and 
B.J.~Pettis' \cite{P} theorem  that a uniformly 
convex space is reflexive. 
While perhaps considered elementary today it illustrated how a geometric 
property can be responsible for a topological property. 
Of course a Banach space can be reflexive without being uniformly convex, 
even under renormings, as shown by M.M.~Day \cite{D2}. 
The problem considered for years by functional analysts was does there 
exist a weaker property of a geometric nature which is equivalent to 
reflexivity. 
In this paper we give an affirmative solution by demonstrating that such a  
property exists. 
The property was suggested in 1971 by V.D.~Milman \cite{Mi} (see also 
\cite{DG2}, problem~IV, p.177). 
We prove that a separable Banach space $X$ is reflexive (if and) only if 
there exists an equivalent norm $\dotnormm$ on $X$ so that 
$$\eqalign{
&\text{whenever a sequence $(x_n)\subseteq X$ satisfies }\cr
&\lim_n \lim_m \normm x_n +x_m\normm = 2\lim_n \normm x_n\normm\cr
&\text{then $(x_n)$ must converge in norm.}\cr}
\tag $*$ $$
The ``if'' part of the characterization follows easily from James' famous 
characterization of reflexivity in terms of the sup of linear 
functionals \cite{J1}. 
Indeed given  $x^*\in X^*$ with $\normm x^*\normm =1$ choose $(x_n)\subseteq 
X$ with $x^*(x_n) \to 1$ and $\normm x_n\normm =1$ for all $n$. 
Then $\lim_m \lim_n \normm x_n + x_m\normm =2$ and so $x_n \to x$ with 
$\normm x\normm =1$. 
Thus $x^* (x) =1$ so $x^*$ attains its norm. 
Hence by \cite{J1} $X$ is reflexive. 

%%% B
The investigation of spaces having property $(*)$ 
 (also called property $(2R)$ in \cite{D1}) 
goes back to the 1950's. 
In \cite{FG}, for example, the relation of $(*)$ to other smoothness and 
rotundity properties was studied. 
For a more complete survey of these notions we refer the reader to \cite{DGZ}. 

More recently (over the past 30 years) functional analysts have considered 
the question as to what sort of nice infinite dimensional subspaces one can 
find in an arbitrary infinite dimensional Banach space $X$. 
One can assume $X$ has a basis and ask what kinds of block bases it has? 
Must one be unconditional? 
Is some block subspace either reflexive or isomorphic to $c_0$ or $\ell_1$? 
These problems are related. 
James \cite{J3} showed that if $(x_i)$ is an unconditional basis for $X$ 
then  either $X$ is reflexive or some block basis is equivalent to the unit 
vector basis for $c_0$ or $\ell_1$. 
W.T.~Gowers \cite{G1} proved the following remarkable dichotomy theorem: 
$X$ contains a subspace $Y$ which either has an unconditional basis or 
is H.I. (hereditarily indecomposable; i.e., if $Z\subseteq Y$ and $Z= V
\oplus W$ then $V$ or $W$ must be finite dimensional). 
Gowers and Maurey \cite{GM} proved that both alternatives are possible. 
Then Gowers \cite{G2} proved that a space need not contain $c_0$, $\ell_1$ or 
a reflexive subspace. 

The search for an answer to this last problem led to much research into 
both characterizations of reflexivity and to the characterization as to when 
$X$ contains isomorphs of $c_0$ or $\ell_1$ (e.g., \cite{J1,2,3},  
\cite{R1,2,3}, \cite{BP}, \cite{M}). 
The proof of our characterization of reflexivity led to additional 
characterizations as to when $X$ contains $c_0$ or $\ell_1$ in terms of 
the asymptotic behavior of sequences in $X$. 

There are two main notions of asymptotic properties in Banach spaces. 
The first is that of a spreading model. 
If $(x_n)$ is bounded in $X$ then by using Ramsey theory (cf.\ \cite{B}, 
\cite{BS}, \cite{O}, \cite{BL}) 
one can extract a subsequence 
$(y_n)$ so that for all $k$ and 
$(a_i)_1^k\subseteq \real$, we have the existence of the 
iterated limit  
$$\lim_{n_1\to\infty} \ldots \lim_{n_k\to\infty} 
\Big\| \sum_{i=1}^k a_iy_{n_i} \Big\| \equiv f(a_1,\ldots,a_k)\ .$$ 
If $(y_n)$ does not converge in the norm topology $f(\cdot)$ is a norm on $c_{00}$, the linear 
space of all finitely supported real valued sequences. 
Let $(e_i)$ be the unit vector basis of $c_{00}$. 
If $(y_n)$ does not converge weakly to a nonzero element of $X$ 
then $(e_i)$ is a basis for $E =[(e_i)]$, the 
completion of $c_{00}$ under $f(\cdot)$. 
In this case we call $(e_i)$ or $E$
the {\it spreading model\/} of $(y_n)$. 
If $(x_i)$ is weakly null then the spreading model $(e_i)$ is unconditional. 
In any event the spreading model is subsymmetric $(\|\sum a_ie_i\| = 
\|\sum a_i e_{n_i}\|$ if $(a_i) \subseteq \real$ and $n_1<n_2 <\cdots$) 
and $(e_1-e_2,e_3-e_4,\ldots)$ is unconditional. 

The second notion of asymptotic structure is due to Maurey, Milman and 
Tomczak-Jaegermann (see \cite{MT}, \cite{MMT}). 
Let $X$ have a basis $(x_i)$. 
For $x,y\in X$ we write $x<y$ if $\max \supp x <\min \supp y$ where if 
$x=\sum a_i x_i$ then $\supp x = \{i:a_i\ne 0\}$. 
$\langle x_i\rangle_{i\in I}$ denotes linear span of $\{x_i: i\in I\}$ and 
$S_{\langle x_i\rangle_{i\in I}}$ denotes the unit sphere of this span. 
Let $n\in \nat$ and let $(w_i)_1^n$ be a normalized basis for some 
$n$ dimensional space. 
We say $(w_i)_1^n \in \{X\}_n$ if 
$$\align 
&\forall\ k_1\in \nat\ \exists\ y_1\in S_{\langle x_i\rangle_{k_1}^\infty} 
\ \forall\ k_2\in \nat\cr 
&\exists\ y_2 \in S_{\langle x_i\rangle_{k_2}^\infty} \cdots \forall\ 
k_n\in \nat\ \exists\ y_n\in S_{\langle x_i\rangle_{k_n}^\infty}
\endalign$$ 
so that $(y_i)_1^n$ is $1+\ep$-equivalent to $(w_i)_1^n$. 
This means that there exist $A,B$ with $AB\le 1+\ep$ so that 
for all $(a_i)_1^n\subset \real$ 
$$A^{-1} \Big\|\sum_1^n a_i y_i\Big\| 
\le \Big\| \sum_1^n a_i w_i\Big\| 
\le B\Big\| \sum_1^n a_i y_i\Big\|\ .$$ 
Note that if $(e_i)$ is a spreading model of a normalized block basis 
of $(x_i)$ then $(e_i)_1^n\in\{X\}_n$ for all $n$. 

Both notions give a more regular structure in general than that possessed by 
the original space $X$. 
They are a joining of the finite and infinite dimensional structures of the 
space. 
Generally only finite dimensional information can be gleaned about $X$ 
from knowledge of its asymptotic structure. 

For example, note that the Schreier 
space $S$ (\cite{CS}, p.1) has a basis having 
a spreading model isometric to $\ell_1$ and yet $S$ is {\it $c_0$ saturated\/} 
(all infinite dimensional subspaces of $S$ contain $c_0$). 
Tsirelson's space $T$ (the dual of Tsirelson's original space \cite{T} as 
described in \cite{FJ}; see also \cite{CS}) has a basis with the property 
that all spreading models are isomorphic to $\ell_1$ and in addition every 
infinite dimensional subspace contains a sequence whose spreading model is 
isometric to $\ell_1$ \cite{OS}. 
Yet $T$ is reflexive. 
We do have the following result which requires a very strong assumption on 
the class of spreading models. 

\proclaim{Theorem \cite{OS}}
If $(x_i)$ is a basis for $X$ and if every spreading model $(e_i)$ of any 
normalized block bases of $(x_i)$ is 1-equivalent to the unit vector basis 
of $\ell_1$ (respectively, $c_0$) then $X$ contains an isomorph of $\ell_1$ 
(respectively, $c_0$). 
\endproclaim 

In this paper we deduce information about the infinite dimensional structure 
of $X$ from knowledge about its asymptotic structure under equivalent norms. 

We shall show that a separable space $(X,\dotNorm)$ can be given a special 
renorming $\dotnormm$ so that certain information about a given spreading 
model $E$ yields information about the infinite dimensional structure of $X$. 
For example if $\normm e_i\pm e_2\normm =2$ (respectively, $\normm e_1+e_2
\normm = 1$) for some spreading model $(e_i)$ 
of a normalized (and respectively, 
weak null) sequence in $X$ then $X$ contains $\ell_1$ (respectively, $c_0$). 
Furthermore we show that a subspace $Y$  of $X$ is reflexive iff 
$Y$ satisfies $(*)$. 

Our main result is the following theorem. 

\proclaim{Main Theorem} 
Every separable Banach space $X$ admits an equivalent strictly convex 
norm $\dotnormm$ with the following properties. 
\roster
\item"a)" If $(x_m) \subseteq X$ is relatively weakly compact and if 
$$\lim_{m\to\infty} \lim_{n\to\infty }\normm x_m+ x_n\normm = 2\lim_{n\to\infty}
\normm x_n\normm$$ 
then $(x_n)$ is norm convergent.
\item"b)" If $(x_n)\subseteq X$ satisfies 
$$\lim_{m\to\infty} \lim_{n\to\infty} \normm x_m\pm x_n\normm = 
2\lim_{n\to\infty} \normm x_n\normm >0$$ 
then some subsequence of $(x_n)$ is equivalent to the unit vector basis 
of $\ell_1$.
\item"c)" If $(x_n)\subseteq X$ is weakly null and satisfies 
$$\lim_{m\to\infty} \lim_{n\to\infty} \normm x_m+x_n\normm 
= \lim_{n\to\infty} \normm x_n\normm >0$$ 
then some subsequence of $(x_n)$ is equivalent to the unit vector basis 
of $c_0$. 
\endroster
\endproclaim 

This theorem is proved in \S 2. 

As a corollary we deduce Milman's suggested characterization of reflexivity. 
In addition we obtain that $X$ contains $\ell_1$ if under all equivalent 
norms, $Y$ admits a normalized basic sequence having a spreading model 
$(e_i)$ satisfying $\normm e_1 \pm e_2\normm =2$. 
In particular if under all equivalent norms $X$ admits a spreading model 
$(e_i)$ which is 1-equivalent to the unit vector basis of $\ell_1$ then $X$ 
contains an isomorph of $\ell_1$. 
If under all equivalent norms $X$ admits a weakly null sequence having 
spreading model $(e_i)$ with $\normm e_1+ e_2\normm =1$ (e.g., if $(e_i)$ 
is 1-equivalent to the unit vector basis of $c_0$) then $X$ contains an 
isomorph of $c_0$. 
From James' proof that $\ell_1$ and $c_0$ are not distortable \cite{J2} one 
obtains that both implications can be reversed. 

In \S3 we present some corollaries discussed briefly in this introduction. 
Our notation is standard as may be found in \cite{LT}. 

\head \S2. Proof of the Main Theorem\endhead 

We first recall the following results of Maurey and Rosenthal.

\proclaim{Theorem (\cite{M},\cite{R1})} 
Let $X$ be a separable Banach space.
\roster
\item"(a)" $X$ is not reflexive if and only if there exists a normalized 
basic sequence $(x_n)\subseteq X$ satisfying for all $x\in X$ and $\alpha,
\beta\ge0$ with $\alpha+\beta=1$, 
$$\lim_{m\to\infty} \lim_{n\to\infty} \|x+\alpha x_n +\beta x_n\| 
= \lim_{m\to\infty} \|x+x_m\|\ .$$
\item"(b)" $X$ contains an isomorph of $\ell_1$ if and only if there exists 
a normalized basic sequence $(x_n)\subseteq X$ such that for all $x\in X$ and 
$\alpha,\beta \in \real$ with $|\alpha| + |\beta|=1$, 
$$\lim_{m\to\infty} \lim_{n\to\infty} \|x+\alpha x_m + \beta x_n\| 
= \lim_{m\to\infty} \|x+x_m\|\ .$$
\item"(c)" $X$ contains an isomorph of $c_0$ iff there exists a normalized 
basic sequence $(x_n)\subseteq X$ such that for all $x\in X$ 
and $\alpha,\beta\in \real$ with $|\alpha|\vee |\beta|=1$, 
$$\lim_{m\to\infty} \lim_{n\to\infty} \|x+\alpha x_m+\beta x_n\| 
= \lim_{m\to\infty} \| x+x_m\|\ .$$
\endroster
\endproclaim

The intuition behind these results and the techniques employed to prove them 
had their origin in \cite{KM} where {\it types\/} were introduced (and 
further developed in \cite{HM}). 
A type $\tau_{(x_n)}$ on $X$ is a function on $X$ defined by a bounded 
sequence $(x_n)\subseteq X$,  
$\tau_{(x_n)}(x) = \lim_{n\to\infty} \|x+x_n\|$. 
Types give information on the asymptotic behavior of a sequence acting on 
the whole space. 
This contrasts with the notion of a spreading model 
which involves only the asymptotic behavior of the sequence $(x_n)$ itself. 
In this paper we characterize the three properties considered in the 
theorem above in terms solely of the asymptotic behavior of the sequences 
themselves. 
The price that must necessarily 
be paid is that we have to consider this behavior under 
all equivalent norms on $X$. 
%%% END OF "F"

%%% START OF "G" 
Let $(X,\dotNorm)$ be a Banach space over $\real$. 
If $x\in X$ we define the {\it symmetrized type\/} norm \newline 
$\dotNorm_x :X\to [0,\infty)$ by 
$$\|y\|_x = \Big\|x\|y\| +y\Big\| + \Big\| x\|y\| -y\Big\| \text{ for } 
y\in X\ .$$ 

\proclaim{Lemma 2.1} 
For all $x\in X$, $\dotNorm_x$ is an equivalent norm on $X$ satisfying 
$2\|y\| \le \|y\|_x \le 2(1+\|x\|) \|y\|$ for all $y\in X$.
\endproclaim 

\demo{Proof} 
The only property not evident is that $\dotNorm_x$ satisfies the triangle 
inequality. 
It is easy to check that for fixed $u,v\in X$ the function $r\mapsto 
\|ru+v\| + \|ru-v\|$ is symmetric and convex  on $\real$ and thus increasing 
on $[0,\infty)$. 
Hence for $y_1,y_2\in X$
$$\eqalignno{
\|y_1 + y_2\|_x  
&= \Big\| x\|y_1 + y_2\| + y_1 +y_2\Big\| 
+ \Big\| x\|y_1 +y_2\| - y_1 -y_2\Big\|\cr 
&\le \Big\| x(\|y_1\| + \|y_2\|) + y_1 +y_2\Big\| 
+ \Big\| x(\|y_1\| +\|y_2\|) - y_1 - y_2\Big\|\cr 
&\le       \Big\| x\|y_1\| + y_1\Big\| +\Big\| x\|y_2\| +y_2\Big\|
+ \Big\| x\|y_1\|-y_1\Big\| +\Big\| x\|y_2\|-y_2\Big\|\cr 
&= \|y_1\|_x + \|y_2\|_x\ .&\qed\cr}$$
\enddemo 

Let $X$ be a separable Banach space. 
It is well known that $X$ admits an equivalent {\it strictly convex\/} 
norm $\dotNorm$, i.e., $\|x\| = \|y\|=1$ and $\|x+y\| =2$ implies that $x=y$. 
Fix a countable dense subset $C$ in $X$ which is closed under rational 
linear combinations. 
Choose $(p_c)_{c\in C} \subseteq (0,\infty)$ so that $\sum_{c\in C} p_c 
(1+\|c\|) <\infty$ for some (and thus for any) equivalent norm on $X$. 
If $\dotNorm$ is an equivalent norm on $X$, 
define $\dotnormm :X\to [0,\infty)$ by 
$\normm x\normm = \sum_{c\in C} p_c\|x\|_c$. 
By Lemma~2.1, $\dotnormm$ is an equivalent norm on $X$. 
Since $0\in C$ and since the sum of a strictly convex norm and any other 
equivalent norm is also strictly convex, $\dotnormm$ is strictly convex. 

\remark{Remark} 
We have assumed that $X$ is a real Banach space. 
Similar results in the complex case can be obtained using 
$$\|y\|_x = \int_0^{2\pi} \Big\| \,\|y\|x + e^{-i\theta} y\Big\| \,
d\theta\ .$$
\endremark 

Our goal is to show that $\dotnormm$ satisfies the main theorem  
if $\|\cdot\|$ is strictly convex. 

\proclaim{Lemma 2.2} 
Let $\dotnormm = \sum_{c\in C} p_c\dotNorm_c$ 
and let $(x_n)\subseteq X$ be $\dotNorm$-normalized. 
\roster
\item"a)" \qquad If $\ds  \lim_{m\to\infty} \lim_{n\to\infty} 
\normm x_m+x_n\normm = 2\lim_{m\to\infty} \normm x_m\normm$
\endroster
then there exists a subsequence $(x'_n)$ of $(x_n)$ satisfying for all 
$y\in X$ and $\beta_1,\beta_2\ge0$  that 
$$\lim_{m\to\infty} \lim_{n\to\infty} \| y+\beta_1 x'_m +\beta_2x'_n\| 
= \lim_{m\to\infty} \|y+ (\beta_1+\beta_2) x'_m\|\ .$$ 
\roster
\item"b)" \qquad If $\ds \lim_{m\to\infty}\lim_{n\to\infty} \normm x_m\pm 
x_n\normm = 2\lim_{m\to\infty} \normm x_m\normm$
\endroster
then there exists a subsequence $(x'_n)$ of $(x_n)$ satisfying for all 
$y\in X$ and $\beta_1,\beta_2\in\real$ with $|\beta_1| +|\beta_2|\ne0$ that 
$$\align 
&\lim_{m\to\infty} \lim_{n\to\infty} \|y+\beta_1 x'_m + \beta_2x'_n\|\cr 
&\qquad = \lim_{m\to\infty} 
\left( \Big\| y{|\beta_1|\over |\beta_1| + |\beta_2|} +\beta_1 x'_m\Big\| 
+ \Big\| y {|\beta_2|\over |\beta_1| +|\beta_2|} +\beta_2 x'_m\Big\|\right)\ .
\endalign$$
\endproclaim 

\demo{Proof} 
a) We may choose $(x'_n)\subseteq (x_n)$ so that for all $c\in C$, $y\in C$ 
and $\beta_1,\beta_2 \in [0,\infty)\cap \que$   the limits 
$$\lim_{m\to\infty} \lim_{n\to\infty} \|y+\beta_1 x'_m + \beta_2 x'_n\|_c$$ 
exist. 
Indeed this is easily done for fixed parameters and then 
one applies a diagonal argument. 

Our hypothesis is that 
$$\lim_{m\to\infty} \lim_{n\to\infty} \sum_{c\in C} p_c \|x'_m+x'_n\|_c 
= 2\lim_m \sum_{c\in C} p_c \|x'_m\|_c \ .$$ 
Since $\|x'_m +x'_n\|_c \le \|x'_m\|_c + \|x'_n\|_c$, 
$\lim_{m\to\infty} \|x'_m\|_c$ exists, 
$\lim_{m\to\infty}\lim_{n\to\infty}  \|x'_m+x'_n\|_c$ exists 
and $p_c \dotNorm_c \le 2 p_c (1+\|c\|)\dotNorm$ for all $c\in C$ we obtain 
since $\sum_{c\in C} p_c (1+\|c\|)<\infty$ 
that for all $c\in C$, 
$$\lim_{m\to\infty} \lim_{n\to\infty} \|x'_m + x'_n\|_c 
= 2\lim_{m\to\infty}\|x'_m\|_c\ .$$ 
In particular taking $c=0$ we obtain that 
$$\lim_{m\to\infty} \lim_{n\to\infty} \|x'_m + x'_n\| =2\ .$$ 
It follows that since $\|x'_n\|=1$, 
$$\lim_{m\to\infty} \lim_{n\to\infty} \|\beta_1 x'_m +\beta_2x'_n\| 
= \beta_1 +\beta_2 
\tag 1$$ 
for all $\beta_1,\beta_2\in [0,\infty)$. 
Similarly we have for all $c\in C$ and $\beta_1,\beta_2\in [0,\infty)$ that 
$$\lim_{m\to\infty} \lim_{n\to\infty} \|\beta_1 x'_m + \beta_2 x'_n\|_c 
= (\beta_1 +\beta_2) \lim_{m\to\infty} \|x'_m\|_c\ .
\tag 2$$ 
Let $y\in C$, $\beta_1,\beta_2 \in [0,\infty) \cap \que$ with 
$\beta_1+\beta_2>0$. 
Setting $c = \frac{y}{\beta_1+\beta_2}$ in (2) we obtain using (1) that 
$$\alignat 2
&\lim_{m\to\infty} \lim_{n\to\infty} \left( \|y+\beta_1 x'_m 
+ \beta_2 x'_n\|  + \|y-\beta_1 x'_m-\beta_2 x'_n\|\right)&\tag 3\cr 
&\qquad =\lim_{m\to\infty} (\beta_1 +\beta_2) 
\left( \Big\| {y\over \beta_1+\beta_2} +x'_m\Big\| 
+ \Big\|{y\over \beta_1 +\beta_2} -x'_m\Big\|\right)\cr 
&\qquad = \lim_m \left( \|y+(\beta_1+\beta_2) x'_m\| +\|y-(\beta_1+\beta_2) 
x'_m\| \right)\ .
\endalignat$$

The triangle inequality yields 
$$\align 
&\lim_{m\to\infty} \lim_{n\to\infty} \|y\pm (\beta_1 x'_m +\beta_2 x'_n)\|\cr 
&\qquad \le \lim_{m\to\infty} \Big\| y{\beta_1\over\beta_1 +\beta_2}\pm 
\beta_1 x'_m\Big\| 
+\lim_{n\to\infty} \Big\| y {\beta_2\over\beta_1+\beta_2}\pm 
\beta_2 x'_n\Big\| \cr 
&\qquad=\lim_{m\to\infty}\beta_1\Big\|{y\over \beta_1+\beta_2}\pm x'_m\Big\|
                       + \beta_2\Big\|{y\over \beta_1+\beta_2}\pm x'_m\Big\|\cr
&\qquad = \lim_{m\to\infty} \| y\pm (\beta_1+\beta_2)x'_m\| 
\endalign$$ 
Thus from (3) we have 
$$\lim_{m\to\infty} \lim_{n\to\infty} \| y+\beta_1 x'_m +\beta_2 x'_n\| 
= \lim_{m\to\infty} \|y +(\beta_1 +\beta_2) x'_m\|\ .$$ 
This proves a) for $y\in C$, $\beta_1,\beta_2 \in \que \cap [0,\infty)$ 
and hence by a density argument we obtain a) in general. 

b) We may assume a) holds for $(x_n)$. 
The only remaining case we need consider is where $\beta_1 >0$ and 
$\beta_2 <0$. 
Actually we shall consider ``$\beta_1 x'_m - \beta_2 x'_n$'' when 
$\beta_1,\beta_2>0$. 
Arguing as above using $\lim_{m\to\infty}\lim_{n\to\infty}\normm x_m-x_n\normm 
= 2\lim_{m\to\infty} \normm x_m\normm$,  
we may assume that $(x_n)$ satisfies 
$$\lim_{m\to\infty} \lim_{n\to\infty} \|\beta_1 x_m - \beta_2 x_n\|_c 
= (\beta_1 +\beta_2) \lim_{m\to\infty} \|x_m\|_c$$ 
for all $c\in C$. 
Letting $y\in C$ and $c= {y\over \beta_1+\beta_2}$ we obtain 
$$\alignat 2
&\lim_{m\to\infty} \lim_{n\to\infty} \left( \| y+\beta_1 x_m -\beta_2 x_n \|
+ \|y-\beta_1 x_m +\beta_2 x_n\|\right) &\tag 4\cr 
&\qquad = (\beta_1 +\beta_2) \lim_{m\to\infty} 
\left(\Big\| {y\over \beta_1+\beta_2} + x_m\Big\| 
+ \Big\| {y\over \beta_1+\beta_2} - x_m\Big\|\right)\cr 
&\qquad = \lim_{m\to\infty} \left( \|y+ (\beta_1+\beta_2) x_m\| 
+ \|y- (\beta_1 +\beta_2) x_m\|\right)\ .
\endalignat $$
Again by the triangle inequality we have 
$$\alignat 2
&\lim_{m\to\infty} \lim_{n\to\infty} \|y\pm(\beta_1 x_m -\beta_2 x_n)\| 
          &\tag 5\\
&\qquad \le \lim_{m\to\infty} \biggl[ \Big\| y {\beta_1\over\beta_1+\beta_2} 
\pm \beta_1x_m\Big\| 
+ \Big\| y{\beta_2\over\beta_1+\beta_2} \mp\beta_2 x_m\Big\| \ .\cr 
\endalignat $$
Since 
$$\alignat 2
 &\lim_{m\to\infty} \biggl[ \Big\| y {\beta_1\over\beta_1+\beta_2} 
+ \beta_1x_m\Big\| 
+ \Big\| y{\beta_2\over\beta_1+\beta_2} - \beta_2 x_m\Big\|  &\tag 6 \\ 
&\qquad\qquad\qquad + \Big\| y{\beta_1\over\beta_1+\beta_2} -\beta_1 x_m\Big\| 
+ \Big\| y{\beta_2\over \beta_1+\beta_2} +\beta_2 x_m\Big\|\biggr] \cr 
&\qquad = \lim_{m\to\infty} \left( \|y+(\beta_1+\beta_2)x_m\| 
+ \|y- (\beta_1 +\beta_2) x_m\|\right)
\endalignat $$
\noindent
it follows from (4) and (5) that we have 
$$\lim_{m\to\infty} \lim_{n\to\infty} \|y + \beta_1x_m -\beta_2 x_n\| 
= \lim_{m\to\infty} \left( \Big\| y{\beta_1\over\beta_1+\beta_2} +\beta_1x_m 
\Big\| + \Big\| y {\beta_2\over\beta_1+\beta_2} -\beta_2 x_m\Big\|\right)$$ 
which completes the proof of b).\qed
\enddemo 

\remark{Remark} 
In the language of types (\cite{R1}, \cite{M}) a) may be restated as 
$$\align
&\text{if }\quad \ds \lim_{m\to\infty} \lim_{n\to\infty} \normm x_m+ x_n
\normm = 2\lim_{m\to\infty} \normm x_m\normm\quad \text{ then}\\
&\text{$(x_m)$ generates an $\ell_1^+$ type on $(X,\dotNorm)$ (equivalently a 
double dual type).}
\endalign$$
The first part of our next lemma is not new (see e.g., \cite{M}, \cite{R1})
but we include the proof for the 
sake of completeness. 
The second part is a slight twist of Maurey's result that a symmetric 
$\ell_1^+$-type yields $\ell_1$. 
In addition the second part of the next lemma establishes that b) of the 
main theorem holds for $\dotnormm$. 
\endremark

\proclaim{Lemma 2.3} 
Let $(x_n)\subseteq X$ be $\dotNorm$-normalized 
and let $\ep_i \subseteq (0,1)$ decrease to $0$. 
\roster
\item"a)" If\quad $\ds \lim_{m\to\infty} \lim_{n\to\infty} \|x+\beta_1 x_m 
+ \beta_2 x_n\| = \lim_{m\to\infty} \|x+(\beta_1+\beta_2)x_n\|$ 
for all $x\in X$ and $\beta_1,\beta_2 >0$
\endroster
then there  exists a subsequence $(x_{n_i})$ of $(x_i)$ satisfying for 
all $1\le i_0\le k$ and $(\alpha_{i_0},\alpha_{i_0+1},\ldots,\alpha_k)
\subseteq [0,\infty)$ that 
$$\Big\| \sum_{i=i_0}^k \alpha_i x_{n_i} \Big\| \ge 
(1-\ep_{i_0}) \sum_{i=i_0}^k \alpha_i\ .$$ 
In particular $(x_n)$ has no weakly null subsequence. 
\roster
\item"b)" If\quad $\ds \lim_{m\to\infty} \lim_{n\to\infty} \|x+\beta_1x_m 
+ \beta_2 x\| = \lim_{m\to\infty} \left( \Big\| y {|\beta_1|\over 
|\beta_1| + |\beta_2|} + \beta_1 x_m\Big\| 
+ \Big\| y{|\beta_2| \over |\beta_1| + |\beta_2|} +\beta_2x_m\Big\|\right)$
for all $x\in X$ and $\beta_1,\beta_2\in\real$ with $|\beta_1| + |\beta_2| 
\ne0$. 
\endroster
then there is a subsequence $(x_{n_i})$ so that for all $1\le i_0\le k$ and  
$(\alpha_i)_{i_0}^k \subseteq \real$, 
$$\Big\| \sum_{i=i_0}^k \alpha_i x_{n_i}\Big\| \ge (1-\ep_{i_0}) 
\sum_{i=i_0}^k |\alpha_i|\ .$$ 
In particular $(x_{n_i})$ is equivalent to the unit vector basis of $\ell_1$. 
\endproclaim 

\demo{Proof} 
a) Given $\delta_i\downarrow 0$ we can 
choose $(x_{n_i})\subseteq (x_i)$ satisfying the following. 
For all $m<\ell$ and $y\in\text{Ball}\langle x_{n_i}\rangle_{i=1}^{m-1}$, 
we have $\|y+\beta_1 x_{n_m} +\beta_2 x_{n_\ell}\|  \ge (1-\delta_m) 
\|y+ (\beta_1+\beta_2) x_{n_m}\|$ if $\beta_1,\beta_2 \in [0,1]$. 
Thus if $(\beta_i)_{i_0}^k \subseteq [0,1]$, $\sum_{i_0}^k \beta_i=1$ 
then 
$$\align
\Big\| \sum_{i=i_0}^k \beta_i x_{n_i}\Big\| 
& = \Big\| \sum_{i=i_0}^{k-2} \beta_i x_{n_i} +\beta_{k-1} x_{n_{k-1}} 
+\beta_2 x_{n_k}\Big\| \cr 
&\ge (1-\delta_{k-1}) \Big\| \sum_{i=i_0}^{k-2} \beta_i x_{n_i} 
+ (\beta_{k-1} +\beta_k) x_{n_{k-1}}\Big\|\cr 
&\ge (1-\delta_{k-1}) (1-\delta_{k-2}) 
\Big\| \sum_1^{k-3} \beta_i x_{n_i} + (\beta_{k-2} +\beta_{k-1} +\beta_k) 
x_{n_{k-2}}\Big\|\cr 
&\ge \cdots \ge \prod_{i=i_0}^{k-1} (1-\delta_i) 
\Big\|\sum_{i=i_0}^k \beta_i x_{n_{i_0}} \big\| 
= \prod_{i=i_0}^{k-1} (1-\delta_i)\ .
\endalign$$
a) follows if we choose the $\delta_i$'s to satisfy $\prod_{i=i_0}^\infty 
(1-\delta_i) \ge 1-\ep_{i_0}$ for all $i_0$. 
The ``in particular'' assertion is immediate from Mazur's theorem. 

b) The argument here is similar but slightly more complicated than a) in 
as much as the condition in b) is not as nice as the one in a). 
Let $\delta_i\downarrow 0$ satisfy $\prod_{i=i_0}^\infty (1-\delta_i) > 
1-\ep_{i_0}$ for all $i_0$ and using the assumption  choose $(x_{n_i}) 
\subseteq (x_i)$ to satisfy for all $m<\ell$ and $y\in \text{Ball}\langle 
x_{n_i}\rangle_{i=1}^{m-1}$, 
$$\|y+\beta_1 x_{n_m}+\beta_2 x_{n_\ell}\| 
> (1-\delta_m) \left[ \Big\| {|\beta_1|\over |\beta_1| + |\beta_2|} 
y+\beta_1 x_{n_m}\Big\| + \Big\| {|\beta_2|\over |\beta_1|+|\beta_2|} y 
+ \beta_2 x_{n_m}\Big\|\right]
\tag 1$$
if $\beta_1,\beta_2\in [-1,1]$ with $|\beta_1| + |\beta_2|\ne0$.

We now show by induction on $k$ that $\|\sum_{i=i_0}^k \beta_ix_{n_i}\| 
\ge \prod_{i=i_0}^{k-1} (1-\delta_i)$ if $i_0\le k$ and $\sum_{i=i_0}^k 
|\beta_i| =1$. 
The claim is trivial for $k=1$ (taking $\prod_\phi (1-\delta_i)\equiv 1$). 
Assume validity of the claim for $k$ and let $\sum_{i=i_0}^{k+1} |\beta_i|=1$. 
For simplicity of the exposition assume $\beta_i\ne0$ for $i_0\le i\le k+1$ 
(the general case follows by a density argument). 
Thus letting $y= \sum_{i=i_0}^{k-1} \beta_i x_{n_i}$ in (1), 
$$\eqalignno{
\Big\|\sum_{i=i_0}^{k+1} \beta_i x_{n_i}\Big\| 
& \ge (1-\delta_k) \biggl[ \Big\| {|\beta_k|\over |\beta_k|+|\beta_{k+1}|} 
\sum_{i=i_0}^{k-1} \beta_i x_{n_i} +\beta_k x_{n_k}\Big\|\cr 
&\qquad\qquad + \Big\| {|\beta_{k+1}| \over |\beta_k| + |\beta_{k+1}|} 
\sum_{i=i_0}^{k-1} \beta_i x_{n_i} +\beta_{k+1} x_{n_k}\Big\|\biggr]\cr 
&\ge (1-\delta_k) \prod_{i=i_0}^{k-1} (1-\delta_i) \biggl[ {|\beta_k| \over 
|\beta_k| +|\beta_{k+1}|} \sum_{i=i_0}^{k-1} |\beta_i| +|\beta_k| \cr 
&\qquad\qquad + {|\beta_{k+1}|\over |\beta_k| +|\beta_{k+1}|} 
\sum_{i=i_0}^{k-1} |\beta_i| + |\beta_{k+1}| \biggr] \cr 
& = \prod_{i=i_0}^k (1-\delta_i)\ . &\qed\cr}$$
\enddemo 

\proclaim{Lemma 2.4} 
Let $\dotnormm = \sum_{c\in C} p_c\|\cdot\|_c$, where $\|\cdot\|$ is an 
equivalent strictly convex norm on $X$.
Let $(x_n)\subseteq X$ be a relatively 
weakly compact sequence. 
If 
$$\lim_{m\to\infty} \lim_{n\to\infty} \normm x_m+x_n\normm 
= 2\lim_{m\to\infty} \normm x_m\normm$$ 
then $(x_n)$ is norm convergent.
\endproclaim 

\demo{Proof} 
Since $\dotnormm$ is a strictly convex norm we need only show that $(x_n)$ 
has a convergent subsequence. 
Indeed if then $(x_n)$ were not convergent it would have two 
subsequences converging to $x\ne y$ respectively. 
But our hypothesis yields $\normm x+y\normm = 2\lim\normm x_m\normm 
= \normm x\normm + \normm y\normm$ which is impossible. 

By passing to a subsequence of $(x_n)$ we may assume that $x_n = x+y_n$ 
where $(y_n)$ is weakly null and $\lim_{n\to\infty} \|y_n\|$ exists. 
If $(y_n)$ were not norm null, we may also assume $\|y_n\| =1$ for all $n$. 
 From Lemma~2.2, passing to a further subsequence, we may assume that for 
all $y\in X$, 
$$\lim_{m\to\infty} \lim_{n\to\infty} \|y+x_m +x_n\| = \lim_{m\to\infty} 
\| y + 2x_m\|\ .$$
For $z\in X$, 
letting  $y= z-2x$ we obtain 
$$\lim_{m\to\infty} \lim_{n\to\infty} \|z+y_m + y_n\| = \lim_{m\to\infty} 
\|z+2y_m\| \ .\tag 1$$ 
Since in particular $\lim_{m\to\infty} \lim_{n\to\infty} \|y_m+y_n\|=2$ it
 follows from (1) and  the definition of $\dotnormm$ that 
$$\lim_{m\to\infty} \lim_{n\to\infty} \normm y_m+y_n\normm 
= 2\lim_{m\to\infty} \normm y_m\normm\ .$$
By Lemma 2.2~a) and Lemma 2.3~a) we conclude that $(y_n)$ is not weakly 
null which is a contradiction.\qed  
\enddemo 

Summarizing our progress thus far we have shown that b) of the main theorem 
is satisfied for $\dotnormm = \sum_{c\in C} p_c \dotNorm_c$  
and in addition a) holds if $\dotNorm$ 
is furthermore a strictly convex norm on $X$. 

\proclaim{Lemma 2.5} 
Let $\dotnormm = \sum_{c\in C} \dotNorm_c$. 
If $(x_n)\subseteq X$ is weakly null and satisfies 
$$\lim_{m\to\infty} \lim_{n\to\infty} \normm x_m + x_n\normm = 
\lim_{m\to\infty} \normm x_m\normm >0$$ 
then $(x_n)$ admits a subsequence which is equivalent to the unit vector 
basis of $c_0$. 
\endproclaim 

\demo{Proof} 
Let $(x_n)$ satisfy the hypothesis of the lemma for $\dotnormm = 
\sum_{c\in C} p_c\dotNorm_c$.  
We may assume $(x_n)$ is basic, $\|x_n\|=1$ for all $n$, and that for all 
$y\in X$ and $\beta_1,\beta_2\in\real$ the following limits exist: 
$$\lim_{m\to\infty} \lim_{n\to\infty} \|y+\beta_1 x_m +\beta_2 x_n\|\ .$$
Since $(x_n)$ is weakly null for all $y\in X$, 
$$\lim_{m\to\infty} \lim_{n\to\infty} \|x_m+ x_n\|_y \ge 
\lim_{m\to\infty} \|x_m\|_y\ .$$ 
As in the proof of Lemma 2.2 since 
$$\lim_{m\to\infty} \lim_{n\to\infty} \sum_{c\in C} p_c \|x_m+ x_n\|_c 
= \lim_{m\to\infty} \sum_{c\in C} p_c \|x_m\|_c$$ 
we obtain for all $y\in C$ and hence in $X$ that 
$$\lim_{m\to\infty} \lim_{n\to\infty} \|x_m+ x_n\|_y 
= \lim_{m\to\infty} \|x_m\|_y\ .
\tag 1$$ 
In particular, $\lim_{m\to\infty} \lim_{n\to\infty} \|x_m+ x_n\|=1$. 
Thus by (1) for all $y\in X$, 
$$\lim_{m\to\infty} \lim_{n\to\infty} (\|y+x_m+x_n\| + \|-y+x_m + x_n\|) 
= \lim_{m\to\infty} (\|y+ x_m\| + \| -y+x_m\|)\ .
\tag 2$$ 
Since 
$$\lim_{m\to\infty} \lim_{n\to\infty} \|\pm y +x_m +x_n\| 
\ge \lim_{m\to\infty} \|\pm y+x_m\|$$
we have from (2) that for all $y\in X$, 
$$\lim_{m\to\infty} \lim_{n\to\infty} \|y+x_m +x_n\| 
= \lim_{m\to\infty} \|y+ x_m\|\ .
\tag 3$$ 

Choose $\ep_i\downarrow0$ with $\prod_1^\infty (1+\ep_i) <2$ and choose, 
using (3), a subsequence $(x_{n_i})$ of $(x_n)$ so that for any integer 
$k\ge0$, $k<i<j$, and $F\subseteq \{1,\ldots,k\}$ then 
$$\Big\| \sum_{s\in F} x_{n_s} + x_{n_i} + x_{n_j}\Big\| 
\le (1+\ep_i) \Big\| \sum_{s\in F} x_{n_s} + x_{n_i}\Big\|\ .$$ 
It follows by iterating this inequality that for all finite $F\subseteq \nat$, 
$\|\sum_{s\in F} x_{n_s}\| \le \prod_1^\infty (1+\ep_i) <2$. 
This implies that $(x_{n_i})$ is equivalent to the unit vector basis 
of $c_0$.\qed
\enddemo 

\remark{Remark} 
The proof yields that for any $\ep>0$ by judiciously choosing the $p_c$'s 
and the original strictly convex norm 
one can choose the norm $\dotnormm$ satisfying the conclusion of the main 
theorem to satisfy for all $x\in X$, 
$$\|x\| \le \normm  x\normm \le (1+\ep)\|x\|\ .$$
\endremark 

We give one final corollary of Lemma 2.5. 
Recall that the summing basis $(s_n)$ for $c_0$ is defined by for all $n$ 
by $s_n = \sum_{i=1}^n e_i$. 

\proclaim{Corollary 2.6} 
Let $\dotnormm = \sum_{c\in C} p_c\dotNorm_c$ 
and let $(x_n) \subseteq X$ satisfy 
$$\align 
&\lim_{n_1\to\infty} \lim_{n_2\to\infty} 
\lim_{n_3\to\infty} \lim_{n_4\to\infty} 
\normm x_{n_1} - x_{n_2} + x_{n_3} - x_{n_4}\normm\cr 
&\qquad = \lim_{n_1\to\infty} \lim_{n_2\to\infty} 
\normm x_{n_1} -x_{n_2}\normm >0\ .
\endalign$$ 
{\rm a)} If $(x_n)$ is weak Cauchy but not weakly convergent then some 
subsequence of $(x_n)$ is equivalent to the summing basis.

\noindent
{\rm b)} If $(x_n)$ is weakly null then some subsequence of $(x_n)$ is 
equivalent to the unit vector basis of $c_0$.
\endproclaim 

\demo{Proof} 
Lemma 2.5 yields the following. 
There exists $C<\infty$ so that for all subsequences of $(x_n)$ there exists 
a further subsequence  $(y_n)$ so that for all finite $F\subseteq \nat$, 
$$\Big\| \sum_{n\in F} (y_{2n} - y_{2n-1}) \Big\| \le C\ .$$
Let 
$$\A = \biggl\{(n_i) \in [\nat]: \text{ for all finite }F\subseteq \nat\ ,\quad 
\Big\| \sum_{i\in F} (x_{n_{2i}} - x_{n_{2i-1}}) \Big\| \le C\biggr\}\ .$$
Here $[\nat]$ denotes the set of all subsequences of $\nat$. 
$\A$ is a Ramsey set (see  e.g., \cite{O}) and thus by our remark above 
there exists $M\in [\nat]$ so that $[M]\subseteq \A$. 
Thus by passing to a subsequence we may assume that if $n_1<\cdots <n_{2k}$ 
then $\|\sum_1^k (x_{n_{2i}} -x_{n_{2i-1}})\|\le C$. 

a) By passing to a subsequence of $x_n$ we may assume that $(x_n)$ is 
basic and moreover $(x_1,x_2-x_1,x_3-x_2,\ldots)$ is seminormalized basic 
(see e.g. \cite{Be, Theorem 8} or \cite{R2}). 
Calling this sequence $(y_n)$ we have that $\|\sum_{n\in F} y_n\| \le 2C 
+ \|x_{n_1}\|$ for all finite $F$, and so $(y_n)$ is equivalent to the unit 
vector basis of $c_0$. 
Hence $(x_n)$ is equivalent to the summing basis: $x_n = \sum_{i=1}^n y_i$. 

b) By Elton's theorem (see \cite{O}) we have that either a subsequence of 
$(x_n)$ is equivalent to the unit vector basis of $c_0$ or some 
subsequence $(y_n)$ of $(x_n)$ satisfies 
$$\lim_{k\to\infty} \Big\| \sum_{i=1}^k (-1)^i y_{n_i}\Big\| =\infty
\ \text{ for all }\ n_1<n_2<\cdots\ .$$ 
From our above remarks we have that a subsequence is the unit vector 
basis of $c_0$.\qed
\enddemo 

\remark{Remark} 
1. We do not know if $X$ can be given a norm $\dotnormm$ satisfying: 
$$\text{if }\quad \lim_{m\to\infty} \lim_{n\to\infty} \normm x_m\pm x_n
\normm = \lim_{m\to\infty} \normm x_m\normm >0$$ 
then some subsequence of $(x_n)$ is equivalent to the unit vector basis 
of $c_0$. 

We can show that this is the case for $\dotnormm $ provided 
in addition one has $\lim_{m\to\infty} \lim_{n\to\infty} \|x_m\pm x_n\| 
= \lim_{m\to\infty} \|x_m\|$. 

However the hypothesis of Corollary~2.6~a) does require the assertion 
that $(x_n)$ not be weakly convergent. 
Indeed if $X$ contains $c_0$ then there exists a normalized sequence 
$(y_n)\subseteq X$ 
which is asymptotically 1-equivalent to the unit vector basis of $c_0$ 
and hence 
$$1=\lim_{n_1\to\infty} \lim_{n_2\to\infty} \|y_{n_1}-y_{n_2}\| = 
\lim_{n_1\to\infty} \lim_{n_2\to\infty} \lim_{n_3\to\infty} 
\lim_{n_4\to\infty} \|y_{n_1}-y_{n_2} + y_{n_3} -y_{n_4}\|\ .$$
Thus $x_n=y+y_n$ satisfies the same condition for any $y\ne 0$ but $(x_n)$ 
admits no basic subsequence. 

2. As we have noted parts b)  and c) of the Main Theorem hold for 
any equivalent norm $\dotNorm$ on $X$ where $\dotnormm 
= \sum_{c\in C} p_c \dotNorm_c$. 
From the proof of Lemma~2.4 it follows that whenever $(x_n)\subseteq X$ 
is relatively weakly compact and satisfies 
$$\lim_{m\to\infty} \lim_{n\to\infty}   \normm x_m+x_n\normm = 
\lim_{m\to\infty} \normm x_m\normm$$ 
then $(x_n)$ is relatively norm compact. 
\endremark 

\head \S3. Corollaries\endhead 

%%% START OF "H"
We now give some corollaries. 
Part a) of Corollary~3.1 yields a positive answer to Milman's 
problem mentioned above. 

\proclaim{Corollary 3.1} 
Let $X$ be a separable Banach space. 
$X$ is reflexive (if and) only if 
there exists an equivalent norm $\dotnormm$ on $X$ 
satisfying for any bounded $(x_n)\subseteq X$ 
\roster
\item"a)" If $\lim_{m\to\infty}\lim_{n\to\infty} \normm x_m+x_n\normm 
= 2\lim_n \normm x_n\normm$ then $(x_n)$ is norm convergent. 
\endroster
Furthermore the norm $\dotnormm$ in {\rm a)} satisfies 
\roster
\item"b)" if $(x_n)$ is weakly null but not norm null then 
$$\lim_{m\to\infty} \lim_{n\to\infty} \normm x_m+x_n\normm > 
\lim_{m\to\infty} \normm x_m\normm$$
provided both limits exist.
\endroster
\endproclaim 

\demo{Proof}
The main theorem ($\,$a), c)$\,$)  yields such a norm if $X$ is reflexive. 
Conversely if a) holds let $x^* \in X^*$ with 
$\normm x^*\normm =1$. 
Choose $(x_n)\subseteq X$, $\normm x_n\normm =1$ with $\lim_{n\to\infty} 
x^* (x_n) =1$. 
It follows that 
$\lim_{m\to\infty} \lim_{n\to\infty} \normm x_m+x_n\normm =2$ 
and so by a), $(x_n)$ converges to some $x$ with $\normm x\normm =1$ and 
$x^* (x)=1$. 
Thus $x^*$ achieves its norm.  
By James' theorem \cite{J1} $X$ must be reflexive.\qed
\enddemo 

\proclaim{Corollary 3.2} 
Let $X$ be a separable Banach space. 
Then there exists an equivalent norm $\dotnormm$ on $X$ such that if $Y$ 
is a subspace  of $X$ then $Y$ is reflexive iff {\rm a)} (and {\rm b)}) 
of Corollary~3.1 hold for all bounded $(x_n)\subseteq Y$. 
\endproclaim 

From b) and c) of the main theorem we obtain 

\proclaim{Corollary 3.3} 
Let $X$ be a separable Banach space. 
The following are equivalent. 
\roster
\item"1)" $X$ contains an isomorph of $\ell_1$ (respectively, $c_0$). 
\item"2)" For all equivalent norms $\dotnormm$ on $X$ there exists a 
normalized sequence in $X$ having 
spreading  model $(e_n)$ which is 1-equivalent to the unit vector basis 
of $\ell_1$ (respectively, $c_0$). 
\item"3)" For all equivalent norms $\dotnormm$ on $X$ there exists a 
normalized (and respectively, weakly null) sequence in $X$ 
having spreading model $(e_n)$ satisfying $\normm e_1 \pm e_2\normm =2$ 
(respectively $\normm e_1 + e_2\normm =1$). 
\endroster
\endproclaim 

In addition to the main theorem the proof requires James' proof that $\ell_1$ 
and $c_0$ are not distortable (\cite{J2} or \cite{LT, p.97}). 
Indeed $1) \Rightarrow 2)$ or 3) is well known from James' result. 
Our discovery is the reverse implications. 

Our work also yields the following corollaries. 

\proclaim{Corollary 3.4} 
Let $X$ be a separable Banach space. 
The following are equivalent. 
\roster
\item"(1)" $X$ is not reflexive.
\item"(2)" For all equivalent norms $\dotnormm$ on $X$ there exists a 
$\dotnormm$ normalized basic sequence $(x_i)$ having spreading model 
$((e_i),\dotnormm)$ satisfying for all $(a_i)\subseteq [0,\infty)$, 
$$\normm \sum a_i e_i\normm = \sum a_i\ .$$
\item"(3)" For all equivalent norms $\dotnormm$ on $X$ there exists a 
$\dotnormm$ normalized basic sequence $(x_i)$ having spreading model 
$((e_i),\dotnormm)$ satisfying 
$$\normm e_1 + e_2\normm = 2\ .$$
\endroster
\endproclaim 

\proclaim{Corollary 3.5} 
Let $X$ be a separable Banach space. 
The following are equivalent.
\roster
\item"(a)" $X$ is reflexive.
\item"(b)" There exists an equivalent norm $\dotnormm$ on $X$ such that 
if $((e_i),\dotnormm)$ is a spreading model of any $\dotnormm$ normalized 
basic sequence in $X$ then $1< \normm e_1 +e_2\normm <2$.
\endroster
\endproclaim 
%%%%%% END OF "H"

\enddocument